\input amstex
\documentstyle{amsppt}
\magnification=\magstep1
\hoffset1 true pc
\voffset2 true pc
\hsize36 true pc
\vsize50 true pc

\NoBlackBoxes
\tolerance=2000 

\define\m1{^{-1}}
\def\gp#1{\langle#1\rangle}
\def\rad{ \operatorname{rad \, }}
\def\Y{\operatorname{ \,\,Y\,\,}}

\TagsOnRight
\topmatter
\title 
On  filtered multiplicative bases of group algebras II
\endtitle 
\author
Victor Bovdi
\endauthor
\dedicatory 
Dedicated to Professor {\it P.M.~Gudivok}  on his 65th birthday 
\enddedicatory 
\leftheadtext\nofrills{Victor Bovdi}
\rightheadtext\nofrills
{
On  filtered multiplicative bases of group algebras II
}

\abstract 
We give an explicit list of all
$p$-groups $G$ with a  cyclic subgroup of index $p^2$, such
that the group algebra $KG$ over the  field $K$ of
characteristic $p$ has a filtered multiplicative $K$-basis. 
We also proved that such  a $K$-basis does not exist  for
the group algebra $KG$, in the case when $G$ is either a
powerful $p$-group or a  two generated $p$-group (
$p\not=2$ ) with  a central cyclic commutator subgroup. This
paper is a continuation  of the related \cite{2}.
\endabstract

\subjclass 
Primary 16A46, 16A26, 20C05. Secondary 19A22 
\endsubjclass

\thanks 
The research was  supported  by OTKA   No.T 025029 and
No.T 029132
\endthanks 

\address
\hskip-\parindent
Institute of Mathematics and Informatics
\newline
University of Debrecen 
\newline
H-4010 Debrecen, P.O. Box 12
\newline
Hungary
\newline
vbovdi\@math.klte.hu
\endaddress

\endtopmatter

\document
\subhead
1. Introduction
\endsubhead
Let $A$ be a finite-dimensional algebra over a field $K$ and let $B$ be a
$K$-basis of $A$.  Suppose that $B$ has the following properties:
\itemitem{1.} if $b_1,b_2\in B$ then either $b_1b_2=0$ or
$b_1b_2\in B$;
\itemitem{2.}  $B\cap \rad(A)$ is a $K$-basis for $\rad(A)$, where 
$\rad(A)$ denotes  the Jacobson radical of
$A$.

Then $B$ is called a {\it  filtered multiplicative
$K$-basis} of $A$.

The filtered multiplicative $K$-basis   arises  in the
theory of representations of algebras and was first
introduced  by H.~Kupisch \cite{5}.  In \cite{1}
R.~Bautista, P.~Gabriel, A.~Roiter and L.~Salmeron  proved
that if there are only finitely many isomorphism classes of
indecomposable $A$-modules over an algebraically closed
field $K$,  then $A$ has   a filtered multiplicative
$K$-basis. Note  that by Higman's theorem  the group
algebra $KG$ over a   field of characteristic $p$ has only
finitely many isomorphism classes of indecomposable
$KG$-modules if  and only if all the Sylow $p$-subgroups of
$G$ are cyclic.

Here  we study  the following question from \cite{1}: {\it
When does  a filtered multiplicative $K$-basis exist  in
the group algebra $KG$}?

Let $G$ be a finite abelian $p$-group. Then
$G=\gp{a_1}\times \gp{a_2} \times \cdots \times \gp{a_s}$
is the direct product of cyclic groups $\gp{a_i}$ of order
$q_i$,   the set
$$
B=\{(a_1-1)^{n_1}(a_2-1)^{n_2}\cdots (a_s-1)^{n_s}\,\,\,\,
\mid\,\,\,\,   0\leq n_i< q_i\} 
$$ 
is  a filtered multiplicative $K$-basis of the group
algebra $KG$ over the field $K$ of characteristic $p$.

Moreover, if $KG_1$ and $KG_2$ have  filtered
multiplicative $K$-bases, which we denote by  $B_1$ and
$B_2$ respectively, then $B_1\times B_2$  is a filtered
multiplicative $K$-basis of the group algebra $K[G_1\times
G_2]$.

P.~Landrock and G.O.~Michler \cite{6} proved  that the
group algebra of the smallest Janko group over a field of
characteristic $2$ does not have a  filtered multiplicative
$K$-basis.

L.~Paris \cite{9} gave  examples of  group algebras  $KG$,
which have  no filtered multiplicative $K$-bases. He also
showed   that if $K$ is a field of characteristic $2$ and
either a) $G$ is a quaternion group of order $8$ and also
$K$ contains a  primitive cube root of the unity or b) $G$
is a dihedral $2$-group,  then $KG$  has a filtered
multiplicative $K$-basis. We  showed in \cite{2}  that  for
the class of all metacyclic $p$-groups the groups mentioned
in the items a) and b) are exactly those, for whose group
algebras L.Paris in \cite{9} presented examples of
multiplicative $K$-bases.

\subhead
2. Results 
\endsubhead
Our main results are the following:
 
\proclaim
{Theorem 1} 
Let $K$ be a field of characteristic $p$  and 
let $G$ satisfy one of the following conditions: 
\itemitem{1.} $G$ is  a powerful $p$-group;
\itemitem{2.} $p$ is odd, $G$ is a $2$-generated $p$-group
with the central cyclic commutator subgroup.

Then the group algebra $KG$ does not have a  filtered multiplicative
$K$-basis.
\endproclaim

\proclaim
{Theorem 2}
Let $K$ be a field of characteristic $p$ and let $G$ be a
nonabelian $p$-group with a  cyclic subgroup of index $p^2$. Then the
group algebra $KG$ possesses a filtered multiplicative
$K$-basis if and only if $p=2$ and  one of the following
conditions is satisfied:
\itemitem{1.} $G$ is either the  dihedral $2$-group or
$D_{2^m}\times C_2$ or the central product   $D_{8} \Y C_4$
of $D_8$ with $C_4$;
\itemitem{2.} $K$  contains a   primitive cube root of the
unity and $G$ is either $Q_8\times C_2$ or $Q_8$;
\itemitem{3.}  $G$ is one  of the following groups: 
\endproclaim
$$
\split
 G_{5}=\langle \,\,\, a,c,d  \,\,\, \mid \,\,\,
                             & a^{2^{m-2}}=c^2=d^2=1,\\
                             & (d,a)=(d,c)=1,\,\,  (c,a)=d \,\,\,
                               \rangle,\text{  with  } m\geq 4;\\                           
G_{13}=\gp{\,\,\,  a,c,d \,\,\, \mid \,\,\, 
                             & a^{2^{m-2}}=c^2=d^2=1,\,\, \\ 
                             & (d,a)=(d, c)=1,\,\, (c,a)=a^{2}d \,\,\,
                             }, \text{  with  } m\geq 5;\\
G_{14}=\gp{\,\,\,  a,c,d \,\,\, \mid \,\,\,
                     & a^{2^{m-2}}=d^2=1,\,\,  c^2=a^{2^{m-3}},\\ 
                     & (d,a)=(d, c)=1,\,\,  (c,a)=a^{2}d  \,\,\,
                     }, \text{  with  } m\geq 5;\\
\endsplit
$$$$
\split
G_{17}=\gp{\,\,\, a,c,d \,\,\, \mid \,\,\,
                     & a^{2^{m-2}}=d^2=c^2=1,  \\
                     & (d,a)=a^{2^ {m-3}},\,\,\, (d,c)=1,\,\,\,   
                       (c,a)=d   \,\,\,
                     }, \text{  with  } m\geq 5;\\
G_{18}=\gp{\,\,\,  a,c,d \,\,\, \mid  \,\,\,
                     & a^{2^{m-2}}=d^2=1,\,\,  c^2=d, \\
                     & (d,a)=a^{2^ {m-3}},\,\,\,   
                       (c,a)=a^{2}d \,\,\,
                     }, \text{  with  } m\geq 4;\\
G_{22}=\gp{ \,\,\, a,c,d \,\,\, \mid \,\,\,
                     & a^{2^{m-2}}=c^2=d^2=1,\,\, (d,a)=1,  \\
                     & (d,c)=a^{2^ {m-3}},\,\,\, 
                       (c,a)=a^{-2^ {m-4}}d \,\,\,
                     }, \text{  with  } m\geq 6;\\
G_{23}=\gp{\,\,\,  a,c,d  \,\,\, \mid  \,\,\,
                     & a^{2^{m-2}}=c^2=d^2=1,\,\,\, (d,a)=1, \\ 
                     & (d,c)=a^{2^{m-3}},\,\,\,  
                       (c,a)=a^{2-2^ {m-4}}d\,\,\,
                     },  \text{ with  } m\geq 6;\\
G_{24}=\gp{\,\,\,  a,c,d \,\,\, \mid  \,\,\,
                     & a^{2^{m-2}}=c^2=d^2=1, \,\,\, (d,c)=1,  \\
                     & (d,a)=a^{2^ {m-3}},\,\, 
                       (c,a)= a^{2- 2^{m-4}}d \,\,
                     }, \text{  with  } m\geq 6;\\
G_{25}=\gp{\,\,\, a,c,d \,\,\, \mid \,\,\,
                     & a^{2^{m-2}}=d^2=1,\,\,\,  
                       c^2=a^{2^{m-3}},\,\,\, (d, c)=1,  \\
                     & (d,a)=a^{2^ {m-3}},\,\, 
                       (c,a)= a^{2-2^{m-4}}d \,\,\,
                     }, \text{  with  } m\geq 5.
\endsplit
$$

\subhead
3. Preliminary remarks and notations
\endsubhead
Let $B$ be  a filtered multiplicative $K$-basis  in a
finite-dimensional $K$-algebra $A$. In the proof of the  main
result we use the following simple properties  of $B$ (see
\cite{2}):
\itemitem{\,(I)\,\,\,\,} $B\cap rad(A)^n$ is a $K$-basis of 
$rad(A)^n$ for all $n\geq 1$.
\itemitem{(II)\,\,} If $u,v\in B\setminus rad(A)^k$ and 
$u\equiv v \pmod{rad(A)^k}$ then $u=v$. 

Recall  that  the {\it Frattini subalgebra}  $\Phi(A)$ of $A$ is
defined as the intersection of all maximal subalgebras of $A$
if those exist, and as $A$ otherwise.  In  \cite{3} it was
shown that if $A$ is a nilpotent algebra over a field $K$,
then $\Phi(A)=A^2$. It follows  that
\itemitem{(III)\,} If $B$ is a filtered multiplicative
$K$-basis  of $A$ and if $B\setminus\{1\}\subseteq rad(A)$,
then all elements of $B\setminus rad(A)^2$ are generators of
$A$ over  $K$.

\bigskip 
Let $K$ be a field of characteristic $p$ and $G$ be a finite
$p$-group.  For $a,b\in G$ we define $a^b=b\m1ab$ and
$(a,b)=a\m1b\m1ab$.  Let  $Q_{2^n}$, $D_{2^n}$ and  $C_{p^n}$
be the {\it generalized quaternion group},  the {\it dihedral } 
$2$-group of order $2^n$, and the {\it cyclic} group of order $p^n$,
respectively.
 
A $p$-group $G$ is  called  {\it  powerful}, if $p=2$ then $G/G^4$ is abelian, 
or if $p>2$ and  $G/G^p$ is abelian.

The ideal 
$
I_K(G)=\big\{ \,\,\, \sum_{g\in G}\alpha_gg\in KG \,\,\,\, \mid\,\,\,\,  
\sum_{g\in G}{\alpha_g}=0\,\,\,\, \big\}
$ 
is called the {\it augmentation} ideal of $KG$ and 
$$
I_K(G)\supset I^2_K(G)\supset\cdots\supset  
I^{s}_K(G)\supset I^{s+1}_K(G)=0. 
$$

Then the  subgroup
$
{\frak D}_n(G)=\big\{\,\,\, g\in G \,\,\,\, \mid \,\,\,\,  g-1\in
I^n_K(G)\,\,\, \big\}
$
is called the $n$th {\it dimensional subgroup}  of $KG$. 

We define the {\it  Lazard-Jennings series }( see  \cite{4} )
$$
M_1(G)\supseteq M_2(G)\supseteq \cdots \supseteq M_t(G)=1
$$ 
as follows:
Put $M_1(G)=G$ and 
$
M_i(G)=\gp{\,\,\, (M_{i-1}(G),G), M^p_{[{\frac{i}{p}}]}(G))\,\,\, },
$
where
\itemitem{---} $[\frac{i}{p}]$ is the smallest integer not 
less than $\frac{i}{p}$;
\itemitem{---} $\big(M_{i-1}(G), G\big)=
               \gp{\,\,\,  (u,v) \,\,\, \mid \,\,\,  
                u\in M_{i}(G),\,\,\,   v\in \,\,\,  G\,\, 
                }$;
\itemitem{---} $M^p_{i}(G)$ is the  subgroup generated by
$p$-powers of the elements of  $M_{i}(G)$.

We know  that for the finite  $p$-groups $M_{i}(G)={\frak D}_i(G)$ for all $i$.
\bigskip

Let  $\Bbb I=\{i\in \Bbb N\mid {\frak D}_i(G)\not={\frak
D}_{i+1}(G)\}$ and let $p^{d_i}$ ($i\in \Bbb I$) be the order
of the elementary abelian $p$-group
$$
{\frak D}_i(G)/{\frak D}_{i+1}(G)=\otimes_{j=1}^{d_i}
\gp{u_{ij}{\frak D}_{i+1}(G)}.
$$
Any element $g\in G$ can be  written uniquely   as 
$$
g=u_{11}^{\alpha_{11}}u_{12}^{\alpha_{12}}\cdots
u_{1d_1}^{\alpha_{1d_1}} u_{21}^{\alpha_{21}} \cdots
u_{2d_2}^{\alpha_{2d_2}}\ldots u_{i1}^{\alpha_{i1}}
\cdots u_{id_i}^{\alpha_{id_i}}\cdots 
u_{s1}^{\alpha_{s1}}\cdots
u_{sd_s}^{\alpha_{sd_s}},
$$
where $i\in \Bbb I$, $u_{ij}\in {\frak D}_i(G)$, $0\leq
\alpha_{ij}<p$, and $s$ is defined  as above.

Elements of the form $w=\prod_{l\in \Bbb
I}(\prod_{k=1}^{d_l}(u_{lk}-1)^{y_{lk}})\in I_K(G)$, where
indicies of its factors are in lexicographic order and $0\leq
y_{lk}<p$, are called  regular elements. Its elements have
weight $\mu(w)=\sum_{l\in {\Bbb I}}(\sum_{k=1}^{d_l}ly_{lk})$.
By Jennings Theorem ( see \cite{4} ), the regular elements of
weight greater than or equal to $t$ constitute an $K$-basis
for the ideal $I^t_K(G)$.  Since ${\frak D}_2(G)=\Phi(G)$,
where $\Phi(G)$ is the Frattini subgroup of $G$, we have that
the set  $\{u_{11},u_{12},\ldots,u_{1d_1}\}$ is a minimal
generator system of $G$.

Clearly, $I_K(G)$ is the radical of $KG$. Suppose  that
$B_1=\{1,\,\, b_1,\ldots,b_{|G|-1} \}$ is a filtered
multiplicative $K$-basis of $KG$. Then $B=B_1\setminus \{1\}$
is a filtered multiplicative $K$-basis  of $I_K(G)$ and
consists  of $|G|-1$ elements. By Jennings Theorem in
\cite{4} $\{\,\,\,(u_{1j}-1)+I^2_K(G)\,\,\, \mid\,\,\,
j=1,\ldots,d_1\,\,\, \}$ form a $K$-basis of
$I_K(G)/I^2_K(G)$.

Put $n=d_1$ and $B\setminus(B\cap I^2_K(G))=\{b_1,b_2,\ldots,b_n\}$. 
Thus
$$
b_k\equiv \sum_{i=1}^n\alpha_{ki}(u_{1i}-1)\pmod {I^2_K(G)},
$$ 
where 
$\alpha_{ki}\in K$ and $\Delta=det(\alpha_{ki})\not=0$. 

Clearly, $z_{ji}=(u_{1j},u_{1i})\in {\frak D}_2(G)$ and  
$z_{ji}-1\in I^2_K(G)$.   
Using the identity
$$
\split
(y-1)(x-1)&=[(x-1)(y-1)+(x-1)+(y-1)](z-1)\\
           &+(x-1)(y-1)+(z-1),
\endsplit\tag1
$$
where $z=(y,x)$,  we obtain  that 
$$
(u_{1j}-1)(u_{1i}-1)=(u_{1i}-1)(u_{1j}-1)+(z_{ji}-1) 
\pmod{I^3_K(G)}.\tag2
$$

Then
$$
\split
b_kb_s\equiv\sum_{i=1}^n\alpha_{ki}\alpha_{si}(u_{1i}-1)^2
&+
\sum_{i,j=1\atop {i<j}}^n(\alpha_{ki}\alpha_{sj}+\alpha_{kj}
\alpha_{si})(u_{1i}-1)(u_{1j}-1)\\
&
+\sum_{i,j=1\atop i<j}^n\alpha_{kj}\alpha_{si}(z_{ji}-1)
\pmod{I^3_K(G)},
\endsplit
\tag3
$$
where $k,s=1,\ldots,n$. 

Let us compute  the dimension subgroups.  
If $p\not=2$ then: 
$$
{\frak D}_1(G)=G,\,\,\,\,\,\,
{\frak D}_2(G)=\Phi(G),\,\,\,\,\,\,
{\frak D}_3(G)=\gp{\,\, ({\frak D}_2(G),G),\,\, G^p\,\, },\,\,\,
$$$$
{\frak D}_4(G)=\gp{\,\, ({\frak D}_3(G),G),\,\,G^p\,\,},\ldots, 
{\frak D}_p(G)=\gp{\,\, ({\frak D}_{p-1}(G),G),\,\, G^p\,\,},
$$
and if $p=2$, we get the following:  
$$
{\frak D}_1(G)=G,\,\,\,\,\,\,
{\frak D}_2(G)=\Phi(G)=G^2,\,\,\,\,\,\,
{\frak D}_3(G)=\gp{\,\, ({\frak D}_2(G),G),\,\, G^4\,\, }.
$$

Assume  that $G$ is  a powerful $p$-group, i.e. if $p=2$ then
$G'<G^4$, and $G'<G^p$  for $p>2$. Then $z_{ji}\in {\frak
D}_3(G)$ and $z_{ji}-1\in I^3_K(G)$.   By (2) it follows that
$b_kb_s\equiv b_sb_k \pmod{I^3_K(G)}$.

Let  
$b_kb_s\in {I^3_K(G)}$. Since the elements  
$$
\{
\,\,
(u_{1i}-1)^2; \,\,
(u_{1k}-1)(u_{1l}-1);\,\,
(u_{2j}-1)
\,\,
\mid 
\,\,
i,k,l=1,\ldots, n; 
\,\, 
k<l; 
\,\, 
j=1,\ldots,d_2\,\,
\}
$$
have  weight $2$, by Jennings Theorem ( see \cite{4} ), these
elements form a basis of $I^2_K(G)$  modulo $I^3_K(G)$.
Because of  (3)  we have that $\alpha_{ki}\alpha_{si}=0$ and
$\alpha_{ki}\alpha_{sj}+\alpha_{si}\alpha_{kj}=0$, where
$i,j=1,\ldots,n$.  Therefore, all minors of  order two, which
are formed from the $k$ and $s$ lines of the matrix
$(\alpha_{i,j})$, equal  zero. From this follows that the
determinant of the matrix $(\alpha_{i,j})$ is zero, which is
impossible.

Therefore, $b_kb_s,b_sb_k\not\in I^3_K(G)$ and $b_kb_s\equiv
b_sb_k \pmod {I^3_K(G)}$ and by property (II) of the filtered
multiplicative $K$-basis  we conclude  that $b_kb_s=b_sb_k$
(for all $k,s=1,\ldots,n$) and $I_K(G)$ is a commutative
algebra, which is a contradiction.

Let $K$ be a field of characteristic $p$, where $p$ is odd,
and let $G$ be a $2$-generated $p$-group with a   central
cyclic commutator subgroup. Then by Theorem $1$ in \cite{7}
$$
\split
G=\gp{\,\,\,  a,c \,\,\,    \mid \,\,\,
 d=(c,a), \,\,\,  
a^{p^{m_1}}=d^{Rp^{r}},\,\,\, 
&c^{p^{m_2}}=d^{Sp^{s}}, \,\,\,
\\
& d^{p^{m_3}}=1, \,\,\,
 (d,a)=(d,c)=1\,\,\,  },
\endsplit
$$ 
where $m_1,m_2,m_3,R,r,S,s$ are  natural numbers defined  in
\cite{7}.

Clearly, $d=(c,a)\in {\frak D}_2(G)$ and   $d-1\in I^2_K(G)$. Then 
$$
(c-1)(a-1)\equiv (a-1)(c-1)+(d-1) \pmod{I^3_K(G)},
$$$$
(c-1)^2(a-1)\equiv (a-1)(c-1)^2+2(c-1)(d-1) \pmod{I^4_K(G)}.
$$

Put 
$$
\cases
b_1\equiv \alpha_1(a-1) +\alpha_2(c-1)\pmod{I^2_K(G)};\\
b_2\equiv \beta_1 (a-1) +\beta_2 (c-1)\pmod{I^2_K(G)},\\
\endcases\tag4
$$
where $\alpha_i,\beta_i\in K$ and
$\Delta=\alpha_1\beta_2- \alpha_2\beta_1\not=0$.

In the rest of the proof we can assume  that $d-1\not\in
I^3_K(G)$.  In the opposite   case, as we have shown  above, we
have a contradiction.  Thus $d-1$ has weight $2$.

Now let us compute $b_{i_1}b_{i_2}b_{i_3}$  modulo
${I^4_K(G)}$, ($i_k=1,2$). The result  of our  computation
will be written  in a  table, consisting of the coefficients
of the decomposition $b_{i_1}b_{i_2}b_{i_3}$ with respect to the  basis
$$
\{\,\,\, (a-1)^{j_1}(c-1)^{j_2}(d-1)^{j_3}
\mid \,\,\, j_1+j_2+2j_3=3;\,\,
j_1,j_2=0,\ldots,3;\,\, j_3=0,1\,\,\, \}
$$
of the ideal $I^3_K(G)$. We will divide the table into  two
parts (the second part written  below  the first part). The 
coefficients corresponding to the first three basis
elements will be in the first part of the table, while the
next three will be in the second one.

$$
\vbox {{ \offinterlineskip
 \halign { 
$#$ \quad  
&  \vrule   \quad   $#$ \quad  
&  \vrule   \quad   $#$ \quad  
&  \vrule   \quad   $#$ \quad  
&  \vrule   \quad   $#$ \quad  
\cr 
       &  (a-1)^3       &   (c-1)^3     &   (a-1)^2(c-1)     \cr
\noalign {\hrule}
b_1b_2b_1                                     & 
\alpha_1^2\beta_1                             & 
\alpha_2^2\beta_2                             &
\alpha_1^2\beta_2+2\alpha_1\alpha_2\beta_1    &  
\cr
b_1b_2^2                                      &
\alpha_1\beta_1^2                             &
\alpha_2\beta_2^2                             &
\alpha_2\beta_1^2+2\alpha_1\beta_1\beta_2     &
\cr
b_2^2b_1                                      &
\alpha_1\beta_1^2                             &
\alpha_2\beta_2^2                             &
\alpha_2\beta_1^2+2\alpha_1\beta_1\beta_2     &
\cr
b_1^2b_2                                      &
\alpha_1^2\beta_1                             &
\alpha_2^2\beta_2                             &
\alpha_1^2\beta_2+2\alpha_1\alpha_2\beta_1    &
\cr
b_2b_1^2                                      &
 \alpha_1^2\beta_1                            &
\alpha_2^2\beta_2                             &
\alpha_1^2\beta_2+2\alpha_1\alpha_2\beta_1    &
\cr
b_2b_1b_2                                     &
\alpha_1\beta_1^2                             &
\alpha_2\beta_2^2                             &
\alpha_2\beta_1^2+2\alpha_1\beta_1\beta_2     &
\cr
b_1^3                                         &
\alpha_1^3                                    &
\alpha_2^3                                    &
3\alpha_1^2\alpha_2                           &
\cr
b_2^3                                         &
\beta_1^3                                     &
\beta_2^3                                     &
3\beta_1^2\beta_2                             &
                   \cr
\noalign {\hrule}
}}}
$$

$$
\vbox {{ \offinterlineskip
 \halign { 
$#$ \quad  
&  \vrule   \quad   $#$ \quad  
&  \vrule   \quad   $#$ \quad  
&  \vrule   \quad   $#$ \quad  
&  \vrule   \quad   $#$ \quad  
\cr 
       &  (a-1)(c-1)^2  &   (a-1)(d-1)  &   (c-1)(d-1)     \cr
\noalign {\hrule}
b_1b_2b_1                                     & 
\alpha_2^2\beta_1+2\alpha_1\alpha_2\beta_2    &
\alpha_1^2\beta_2                             &  
\alpha_2^2\beta_1+2\alpha_1\alpha_2\beta_2    &
\cr
b_1b_2^2                                      &
\alpha_1\beta_2^2+2\alpha_2\beta_1\beta_2     &
 \alpha_1\beta_1\beta_2                       &
3\alpha_2\beta_1\beta_2                       &
\cr
b_2^2b_1                                      &
\alpha_1\beta_2^2+2\alpha_2\beta_1\beta_2     &
3\alpha_1\beta_1\beta_2                       &
\alpha_2\beta_1\beta_2+2\alpha_1\beta_2^2     &
\cr
b_1^2b_2                                      &
\alpha_2^2\beta_1+2\alpha_1\alpha_2\beta_2    &
3\alpha_1\alpha_2\beta_1                      &    
2\alpha_2^2\beta_1+\alpha_1\alpha_2\beta_2    &
\cr
b_2b_1^2                                      &
 \alpha_2^2\beta_1+2\alpha_1\alpha_2\beta_2   &
\alpha_1\alpha_2\beta_1                       &    
3\alpha_1\alpha_2\beta_2                      &
\cr
b_2b_1b_2                                     &
\alpha_1\beta_2^2+2\alpha_2\beta_1\beta_2     &
\alpha_2\beta_1^2                             &    
\alpha_1\beta_2^2+2\alpha_2\beta_1\beta_2     &
\cr
b_1^3                                         &
3\alpha_1\alpha_2^2                           &
3\alpha_1^2\alpha_2                           &    
3\alpha_1\alpha_2^2                           &
\cr
b_2^3                                         &
3\beta_1\beta_2^2                             &
3\beta_1^2\beta_2                             &    
3\beta_1\beta_2^2                             &
                   \cr
\noalign {\hrule}
}}}
$$
We have obtained $8$ elements.  If $char K>3$ or $char K=3$
and $|a|\not=3\not=|c|$,  then the $K$-dimension of
$I^3_K(G)/I^4_K(G)$ equals $6$.

In the opposite case the  $K$-dimension of
$I^3_K(G)/I^4_K(G)$ is less than $6$. But we must have
either $6$ or less, respectively, linearly independent
elements $b_{i_1}b_{i_2}b_{i_3}$  modulo the ideal
$I^4_K(G)$.  Therefore, some of these elements either
are  equal to zero modulo the ideal $I^4_K(G)$ or
coincide with some other elements of the system. Then we
get by simple calculations that $\Delta=0$, which is
impossible  for $charK>2$.

Indeed,  for example, if  $b_1b_2b_1\equiv 0\pmod{I^4_K(G)}$  
then  from the first two columns follows that either 
$\alpha_2=\beta_1=0$ or  $\alpha_1=\beta_2=0$. Also,   
from  the 3th and  6th columns we get either  
 $\alpha_1^2\beta_2=0$ or $\alpha_2^2\beta_1=0$, respectively. 
Therefore,  $\Delta=0$, which is impossible.

Now, for example, put   
$\beta_1\beta_2\beta_1\equiv\beta_1\beta_2^2\pmod{I^4_K(G)}$. 
It follows that     
$\alpha_1=\beta_1$, $\alpha_2=\beta_2$ and $\Delta=0$, which 
is also a contradiction. 

The rest of the cases are similar to these two.

\subhead
4. Proof of Theorem 2
\endsubhead
According to  Theorems $1$ and $2$ in  \cite{8} there are
$11$ finite nonabelian $p$-groups  of order $p^m$ and
exponent $p^{m-2}$ for  $p>2$ and $26$ such groups when
$p=2$.

In the rest of the proof we will keep the indexes of these
groups as  in \cite{8}, but to make the calculation easier
we will rename the generators  of these groups as
follows: we do not change  $a$ and $c$ but   we will denote
$d$ either as  $b$ or $b\m1$.

First, we consider the case when  $charK>2$. Then  by
Theorem $1$ in  \cite{8} a finite nonabelian neither
metacyclic  nor powerful $p$-group of order $p^m$ and
exponent $p^{m-2}$ is isomorphic to one of the groups
$G_{1}$, $G_{5}$, $G_{6}$, $G_{7}$, $G_{11}$.
$$
\split
\text{ If  } G=G_{1}=\gp{\,\,\, a,c,d \,\,\, 
                    \mid \,\,\,  
                    & a^{p^{m-2}}=c^p=d^p=1, \\ 
                    & (d,a)=1,  \,\,\,
                      (d,c)=1,\,\,\, 
                      (c,a)=d \,\,\,
            },   
\text{  with }  m\geq 3,
\endsplit
$$
then by Theorem $1$ we have  that $KG$ has no  filtered
multiplicative $K$-basis.

Let $G$ be  one of the following groups:
$$
\split
H(r)=\gp{\,\,\, a,c,d  \,\,\, 
                      \mid \,\,\,  
                    & a^{p^{m-2}}=c^p=d^p=1,\,\,\,(d,a)=1,\\
                    & (c,a)=d,\,\,\,  
                      (d,c)=a^{-rp^{m-3}} \,\,\, 
                    }, \text{ with  } m\geq 4; \\
G_{7}=\gp{\,\,\, a, c, d \,\,\, 
                      \mid \,\,\,  
                   &  a^{p^{m-2}}=c^p=d^p=1, \\    
                   &  (d,a)=a^{p^{m-3}}, \,\,\,   
                      (c,a)=d,\,\,\,  
                      (d,c)=1\,\,\,  
                   },  \text{ with } m\geq 4;\\
G_{11}=\gp{\,\,\,  a,c,d \,\,\, 
                       \mid \,\,\,  
                     & a^9=d^3=1,\\  
                     & a^3=c^3,\,\,  
                       (d,a)=1,\,\,\,      
                       (c,a)=d,\,\,\,  
                       (d,c)=a^3\,\,\,
                     },\\
\endsplit
$$
where either  $r=1$ or $r$ is a quadratic nonresidue modulo
${p}$.  Note that if  $r=1$ then $H(1)$ coincides with
$G_{5}$ and in the other case $H(r)$ coincides with $G_6$.

Then by (1) we get  
$$
(c-1)(a-1)\equiv (a-1)(c-1)+(d-1) \pmod{I^3_K(G)}.
$$  
Now, similarly to the proof of Theorem $1$  above, let us
compute $b_{i_1}b_{i_2}b_{i_3}$  modulo ${I^4_K(G)}$,
($i_k=1,2$). The result  of our computation  can be written
in a similar table, consisting of the coefficients of the
decomposition $b_{i_1}b_{i_2}b_{i_3}$ with respect to the 
considered basis 
$$
\{\quad 
(a-1)^3,     \quad
(a-1)^2(c-1),\quad
(a-1)(d-1),  \quad
(c-1)(d-1),  \quad
(a-1)(c-1)^2 \quad 
\}
$$ 
of the ideal $I^3_K(G)$.  

We must consider the following tree  cases.

Case 1. Let either $G\not=G_7$ with $p=3$ and $m=4$ or
$G\not=H(r)$ with $p=3$, $m=4$ or  $G\not=G_{11}$. Then by
(1), the element $d-1$ is central  modulo $I^4_K(G)$,
and we will have the same table as  in the proof of Theorem 1.
Thus, we will have a contradiction.

Case 2. Let $G=G_7$ with  $p=3$ and $m=4$. 
Then 
$$
G=\gp{\quad a,c,d \quad \mid \quad a^9=c^3=d^3=1,\quad  
(d,a)=a^3,\quad  (c,a)=d,\quad  (d,c)=1\quad } 
$$ 
and by (1) we get $(d-1)(c-1)\equiv (c-1)(d-1) \pmod{I^4_K(G)}$ and 
$$
 (d-1)(a-1)\equiv (a-1)(d-1)+(a-1)^3 \pmod{I^4_K(G)}.
$$
 
Modulo $I^4_K(G)$ it follows   that 
$$
\eightpoint
\vbox {{ \offinterlineskip
 \halign { 
$#$
&  \vrule  \quad     $#$   
&  \vrule  \quad     $#$   
&  \vrule  \quad     $#$   
&  \vrule  \quad     $#$  
&  \vrule  \quad     $#$   
\cr 
 &(a-1)^3 &(a-1)^2(c-1) &(a-1)(d-1) &(c-1)(d-1) &(a-1)(c-1)^2 \cr
\noalign {\hrule}
b_1b_2b_1                                           & 
\alpha_1\beta_1 (\alpha_1+\alpha_2)                 & 
\alpha_1\Delta                                      &
\alpha_1^2\beta_2                                   &
-\alpha_2\Delta                                     & 
-\alpha_2\Delta           
\cr
b_1b_2^2                                            &
\beta_1^2(\alpha_1+\alpha_2)                        & 
-\beta_1\Delta                                      &
\alpha_1\beta_1\beta_2                              &
0                                                   &
\beta_2\Delta            
\cr
b_2b_1^2                                            &
\alpha_1\beta_1(\beta_1+\beta_2)                    &
\alpha_1\Delta                                      &
\alpha_1\alpha_2\beta_1                             &
0                                                   &
-\alpha_2\Delta           
\cr
b_2b_1b_2                                           &
\alpha_1\beta_1(\beta_1+\beta_2)                    &
\beta_1\Delta                                       &
\alpha_2\beta_1^2                                   &
\beta_2\Delta                                       & 
\beta_2\Delta            
\cr
b_1^3                                               &
\alpha_1^2(\alpha_1+\alpha_2)                       & 
0                                                   &
\alpha_1^2\alpha_2                                  &
0                                                   &
0                        
\cr
b_1^2b_2                                            &
\alpha_1\beta_1(\alpha_1+\alpha_2)                  & 
\alpha_1\Delta                                      &
\alpha_1\alpha_2\beta_1                             &
\alpha_2\Delta                                      & 
-\alpha_2\Delta          
\cr
b_2^2b_1                                            &
\alpha_1\beta_1(\beta_1+\beta_2)                    &
-\beta_1\Delta                                      &
\alpha_1\alpha_1\beta_2                             &
-\beta_2\Delta                                      & 
\beta_2\Delta           
\cr
b_2^3                                               &
\beta_1^2(\beta_1+\beta_2)                          &
0                                                   &
\beta_1^2\beta_2                                    &
0                                                   &
0                                                  
\cr
\noalign {\hrule}
}}}
$$
We have obtained $8$ elements, but the $K$-dimension of
$I^3_K(G)/I^4_K(G)$ equals $5$. Therefore, some of these
elements either are  equal to zero modulo the ideal
$I^4_K(G)$ or coincide with some other elements of the
system.  Then we get by simple calculations that
$\Delta=0$, which is impossible.

Case 3. Let either  $G=H(r)$ with $p=3$, $m=4$ or
$G=G_{11}$.  Then  by (1) we get
$$
\split
 (d-1)(c-1)& \equiv (c-1)(d-1)+(a-1)^3 \pmod{I^4_K(G)},\\
 (d-1)(a-1)& \equiv (a-1)(d-1) \pmod{I^4_K(G)}.\\
\endsplit
$$
Modulo $I^4_K(G)$ it follows    that 
$$
\eightpoint
\vbox {{ \offinterlineskip
 \halign { 
$#$
&  \vrule  \quad     $#$   
&  \vrule  \quad     $#$   
&  \vrule  \quad     $#$   
&  \vrule  \quad     $#$  
&  \vrule  \quad     $#$   
\cr 
 &\scriptstyle{(a-1)^3}  
 &\scriptstyle{(a-1)^2(c-1)}   
 &\scriptstyle{ (a-1)(d-1)}  
 &\scriptstyle{(a-1)(c-1)^2} 
 &\scriptstyle{(c-1)(d-1)}  
 \cr
\noalign {\hrule}
b_1b_2b_1                                           & 
\beta_1 (\alpha_1^2+\alpha_2^2)                     & 
\alpha_1\Delta                                      &
\alpha_1\Delta                                      &
-\alpha_2\Delta                                     & 
-\alpha_2\Delta           
\cr
b_1b_2^2                                            &
\beta_1(\alpha_1\beta_1+\alpha_2\beta_2)            & 
-\beta_1\Delta                                      &
\beta_1\Delta                                       &
\beta_2\Delta                                       &
0            
\cr
b_2b_1^2                                            &
\alpha_1(\alpha_1\beta_1+\alpha_2\beta_2)           &
\alpha_1\Delta                                      &
-\alpha_1\Delta                                     &
-\alpha_2\Delta                                     &
0              
\cr
b_2b_1b_2                                           &
\alpha_1(\beta_1^2+\beta_2^2)                       &
-\beta_1\Delta                                      &
-\beta_1\Delta                                      &
\beta_2\Delta                                       & 
\beta_2\Delta            
\cr
b_1^3                                               &
\alpha_1(\alpha_1^2+\alpha_2^2)                     & 
0                                                   &
0                                                   &
0                                                   &
0                        
\cr
b_1^2b_2                                            &
\alpha_1(\alpha_1\beta_1+\alpha_2\beta_2)           & 
\alpha_1\Delta                                      &
0                                                   &
-\alpha_2\Delta                                     & 
\alpha_2\Delta          
\cr
b_2^2b_1                                            &
\beta_1(\alpha_1\beta_1+\alpha_2\beta_2)            &
-\beta_1\Delta                                      &
0                                                   &
\beta_2\Delta                                       & 
-\beta_2\Delta           
\cr
b_2^3                                               &
\beta_1(\beta_1^2+\beta_2^2)                        &
0                                                   &
0                                                   &
0                                                   &
0                                                  
\cr
\noalign {\hrule}
}}}
$$
Now, similarly to case 2, we have obtained $8$ elements,
but the $K$-dimension of $I^3_K(G)/I^4_K(G)$ equals $5$.
Therefore, some of these elements either identically equal
zero modulo the ideal $I^4_K(G)$ or coincide with some
other elements of the system.  Then we get by simple
calculations that $\Delta=0$, which is impossible.

Therefore, if  $p$ is odd, then $KG$ has no  filtered
multiplicative $K$-basis.

Now let $charK=2$.   Then by Theorem $2$ in  \cite{8} a
finite nonabelian neither  metacyclic  nor powerful
$2$-group  of order $2^m$ and exponent $2^{m-2}$ is
isomorphic to one of the  groups $G_2$, $G_{3}$, $G_5$,
$G_{11}$, $G_{12}$, $G_{13}$, $G_{14}$, $G_{15}$, $G_{16}$,
$G_{17}$, $G_{18}$, $G_{22}$, $G_{23}$, $G_{24}$, $G_{25}$,
$G_{26}$ with the exception of  $G_4$, which is discussed below. 
If $m=5$ then $G_{25}$ coincides with $G_{26}$.

First, let us suppose that $G$ is   
$$
\split
G_4=\gp{\,\,\, a,c,d  \,\,\, \mid \,\,\,
                     & a^{2^{m-2}}=c^2=d^2=1, \\
                     & (d,a)=(c,a)=1,\,\,\, 
                       (d,c)=a^{2^{m-3}} \,\,\,
                     },  \text{ with } m\geq 4.\\ 
\endsplit
$$ 
Clearly,   for $m\ge 5$,  $G$ is  powerful, therefore $G_4$
is the  central product of $D_8=\gp{a,b}$ and $C_4=\gp{c}$.
Put $b_1=(1+a)+(1+c)$, $b_2=1+c$ and $b_3=1+d$. Thus $\{
\,\,\, b_2^{i}\,\, b_1^{j}\,\, b_2^{k}\,\, b_3^l  \,\,\,\,
\mid \,\,\,\, i,k=0,1; \,\,\,\, j,l=0, \ldots, 3 \,\,  \}$
form a filtered multiplicative $K$-basis of $KG_4$.

Let $G$ be  either $G_2=Q_{2^{n-1}}\times C_2$ or
$G_3=D_{2^{m-1}}\times C_2$.  Then using  \cite{2}, 
one obtains that the  group
$G$  satisfies   conditions 1. or 2. of Theorem $2$.

Let $G$ be one of the following groups: 
$$
\split
G_{11}=\gp{\,\,\, a,c,d  \,\,\, \mid \,\,\,
                     & a^{2^{m-1}}=c^2=d^2=1,\,\,\, (d,c)=1, \\
                     & (d,a)=1,\,\,\, 
                       (c,a)=a^{2+2^{m-2}} \,\,\,
                     },  \text{ with } m\geq 4;\\ 
G_{12}=\gp{\,\,\, a,c,d  \,\,\, \mid \,\,\,
                     & a^{2^{m-2}}=c^2=d^2=1, \,\,\, (d,c)=a^{2^{m-3}},\\
                     & (d,a)=1,\,\,\,  
                       (c,a)=a^{2} \,\,\,
                     },  \text{ with } m\geq 5;\\ 
\endsplit
$$$$
\split
G_{15}=\gp{\,\,\, a,c,d  \,\,\, \mid \,\,\,
                     & a^{2^{m-2}}=c^2=d^2=1,\,\,\, (d,c)=1, \\
                     & (d,a)=a^{2^ {m-3}},\,\,\, 
                       (c,a)=a^{2+2^{m-3}}\,\,\,  
                     }, \text{  with  } m\geq 5;\\
G_{16}=\gp{\,\,\, a,c,d  \,\,\, \mid \,\,\,
                     & a^{2^{m-2}}=c^2=d^2=1, \,\,\,
                       (d,c)=a^{2^{m-3}}, \\
                     & (d,a)=a^{2^{m-3}},\,\,\,
                       (c,a)=a^{2+2^{m-3}} \,\,\,
                     }, \text{ with } m\geq 5.\\                  
\endsplit
$$   
It is easy to see that $G$ is $3$-generated and we can put 
$$
\cases
b_1 \equiv \alpha_1(1+a) +
    \alpha_2(1+c)+\alpha_3(1+d)   \pmod{I^2_K(G)};\\
b_2 \equiv \beta_1 (1+a) +
    \beta_2 (1+c)+\beta_3(1+d)    \pmod{I^2_K(G)};\\
b_3 \equiv \gamma_1 (1+a) +
    \gamma_2 (1+c)+\gamma_3(1+d)  \pmod{I^2_K(G)},\\
\endcases
$$
where $\alpha_i,\beta_i,\gamma_i\in K$ and
$\Delta\not=0$.

By (1) we have  that 
$
(1+c)(1+a)\equiv (1+a)(1+c)+(1+a)^2 \pmod{I^3_K(G)}
$
and   $1+d$ is  central by modulo ${I^3_K(G)}$.

Now let us compute $b_{i}b_{j}$ by modulo ${I^3_K(G)}$,
($i,j=1,2$). The result  of our  computation will be
written  in a  table, consisting of the coefficients of the
decomposition $b_{i}b_{j}$ with respect to the  basis
$
\{\,\,\,
(1+a)^2,  (1+a)(1+c), (1+a)(1+d), (1+c)(1+d)
\,\,\, \}
$
of the ideal $I^2_K(G)$. 

$$
\vbox {{ \offinterlineskip
 \halign { 
$#$ \quad  
&  \vrule   \quad   $#$  
&  \vrule   \quad   $#$   
&  \vrule   \quad   $#$   
&  \vrule   \quad   $#$   
\cr 
&  (1+a)^2    &(1+a)(1+c) &  (1+a)(1+d) & (1+c)(1+d) \cr
\noalign {\hrule}
b_1b_2                                        & 
\beta_1(\alpha_1+\alpha_2)                    & 
\alpha_1\beta_2+\alpha_2\beta_1               &
\alpha_1\beta_3+\alpha_3\beta_1               &  
\alpha_2\beta_3+\alpha_3\beta_2                 
\cr
b_2b_1                                        &
\alpha_1(\beta_1+\beta_2)                     & 
\alpha_1\beta_2+\alpha_2\beta_1               &
\alpha_1\beta_3+\alpha_3\beta_1               &  
\alpha_2\beta_3+\alpha_3\beta_2                 
\cr
b_1b_3                                        & 
\gamma_1(\alpha_1+\alpha_2)                   & 
\alpha_1\gamma_2+\alpha_2\gamma_1             &
\alpha_1\gamma_3+\alpha_3\gamma_1             &  
\alpha_2\gamma_3+\alpha_3\gamma_2                 
\cr
b_3b_1                                        & 
\alpha_1(\gamma_1+\gamma_2)                   & 
\alpha_1\gamma_2+\alpha_2\gamma_1             &
\alpha_1\gamma_3+\alpha_3\gamma_1             &  
\alpha_2\gamma_3+\alpha_3\gamma_2                 
\cr
b_2b_3                                        &
\gamma_1(\beta_1+\beta_2)                     & 
\beta_1\gamma_2+\beta_2\gamma_1               &  
\beta_1\gamma_3+\beta_3\gamma_1               &  
\beta_2\gamma_3+\beta_3\gamma_2                 
\cr
b_3b_2                                        &
\beta_1(\gamma_1+\gamma_2)                    &  
\beta_1\gamma_2+\beta_2\gamma_1               &  
\beta_1\gamma_3+\beta_3\gamma_1               &  
\beta_2\gamma_3+\beta_3\gamma_2                 
\cr
b_1^2                                        &
\alpha_1(\alpha_1+\alpha_2)                  &
0                                            &
0                                            &
0
\cr
b_2^2                                         &
\beta_1(\beta_1+\beta_2)                      &
0                                             &
0                                             &
0
\cr
b_3^2                                         &
\gamma_1(\gamma_1+\gamma_2)                   &
0                                             &
0                                             &
0
\cr
\noalign {\hrule}
}}}
$$
We have obtained $9$ elements.  But the $K$-dimension of
$I^2_K(G)/I^3_K(G)$ equals $4$ and we must have  four 
linearly independent elements $b_{i}b_{j}$  modulo  ideal
$I^3_K(G)$.  Therefore, some of these elements either
are  equal to zero modulo the ideal $I^3_K(G)$ or
coincide with some other elements of the system.

We will  consider the following cases:

Case 1. Let $b_ib_j\equiv b_jb_i \pmod{I^3_K(G)}$. 
But  $b_ib_j,b_jb_i\not\in I^3_K(G)$ 
and by property (II)
we conclude  that  $I_K(G)$ is a commutative algebra, which
is a contradiction.

Case 2. It is easy to see that the first six lines equal neither  
 zero nor the last  three lines, because 
from the last three columns in our table follows that 
all minors equal zero and $\Delta=0$, which is impossible.

Case 3. Let $b_i^2\equiv 0\pmod{ I_K^3(G)}$ 
for two values  $i$, for example,
$b_1^2\equiv b_2^2\equiv 0\pmod{ I_K^3(G)}$.  Then either $\alpha_1=0$ or
$\alpha_1=\alpha_2$ and either $\beta_1=0$  or
$\beta_1=\beta_2$, respectively.  
Since we must have $4$ linearly
independent  elements and the cases 1 and 2 are impossible,
we have  that $b_3^2\equiv 0\pmod{ I_K^3(G)}$. Thus 
either $\gamma_1=0$ or $\gamma_1=\gamma_2$. If we put the  
values of $\alpha_j,\beta_j$ and $\gamma_j$ into our table ( we will have 
eight cases)  then 
we will get  a contradiction.

Therefore, in any other case, we can refer back to  these 
three cases. 

Suppose  $G$ is  one of the groups $G_5$, $G_{17}$, $G_{22}$ or  
$G_{25}$ with $m=5$ ( i.e. $G_{26}$ ) from  Theorem $2$.  Clearly,  
$G$ is two-generated and we can
assume that $u=b_1$,  $v=b_2$, where $b_1,b_2$ can be written
as in  (4).

Then by (1) we get $1+d$ is central  modulo  $I^3_K(G)$ and 
$$
(1+c)(1+a)\equiv (1+a)(1+c)+(1+d) \pmod{I^3_K(G)}.
$$   
It follows  that 
$$
\cases
{u}{v} \equiv \alpha_1\beta_1(1+a)^2+\Delta(1+a)(1+c)+\alpha_2\beta_1(1+d)  
&\pmod{I^3_K(G)};\\
{v}{u} \equiv \alpha_1\beta_1(1+a)^2+\Delta(1+a)(1+c)+\alpha_1\beta_2(1+d)
&\pmod{I^3_K(G)};\\
{u}^2\equiv \alpha_1^2(1+a)^2+\alpha_1\alpha_2(1+d)
&\pmod{I^3_K(G)};\\
{v}^2\equiv \beta_1^2 (1+a)^2+\beta_1 \beta_2 (1+d)
&\pmod{I^3_K(G)}.\\
\endcases
$$
Since the $K$-dimension of $I^2_K(G)/I^3_K(G)$ equals $3$ and
$\Delta\not= 0$, we have  that ${u}{v}\not\equiv
{v}{u}\pmod{I^3_K(G)}$ and ${u}^2\not\equiv
{v}^2\pmod{I^3_K(G)}$. Thus either ${v}^2\equiv
0\pmod{I^3_K(G)}$ and $\beta_1=0$ or ${u}^2\equiv
0\pmod{I^3_K(G)}$ and $\alpha_1=0$. It is easy to see  that
the second case is symmetric to the first one, so we consider
only the first one.  Therefore $\alpha_1\not=0$ and we can
put ${u}= (1+a)+\mu(1+c)$ and ${v}= 1+c$, where
$\mu=\frac{\alpha_2}{\alpha_1}$.

We will prove  that the elements $\{{u}^{i}, {v}{u}^{i-1},
{v}{u}{v}{u}^{i-3}, {u}{v}{u}^{i-2}\}$ form a basis of
$I^{i}_K(G)/I^{i+1}_K(G)$, where $i>3$.
 
First of all,  ${v}{u}{v}\equiv(1+c)(1+d)\pmod{I^4_K(G)}$, 
and by induction  we get   that  
$$
\split
{u}^{4i}   & \equiv    (1+a)^{4i}                      
           \pmod{I^{4i+1}_K(G)},\\
{u}^{4i+1} & \equiv    (1+a)^{4i+1}+\mu(1+a)^{4i}(1+c) 
           \pmod{I^{4i+2}_K(G)},\\
{u}^{4i+2} & \equiv    (1+a)^{4i+2}+\mu(1+a)^{4i}(1+d) 
           \pmod{I^{4i+3}_K(G)},\\
{u}^{4i+3} & \equiv    (1+a)^{4i+3}+\mu(1+a)^{4i+2}(1+c)\\
         & +     \mu (1+a)^{4i+1}(1+d)+
\mu^2(1+a)^{4i}(1+c)(1+d)
           \pmod{I^{4i+4}_K(G)}.
\endsplit
$$
We consider the following four  cases:  

Case 1. If  $i\equiv 0 \pmod{4}$, then  we have 
$$
\split
{u}^{i}        &  \equiv(1+a)^{i}
                \pmod{I^{i+1}_K(G)}, \\
{v}{u}^{i-1}     &  \equiv(1+a)^{i-1}(1+c)+
                (1+a)^{i-2}(1+d)\\
             &  +(1+a)^{i-3}(1+c)(1+d) 
                \pmod{I^{i+1}_K(G)},\\
({u}{v}){u}^{i-2 } &  \equiv(1+a)^{i-1}(1+c)+
                \mu(1+a)^{i-3}(1+c)(1+d) 
                \pmod{I^{i+1}_K(G)},\\
({v}{u}{v}){u}^{i-3} &  \equiv(1+a)^{i-3}(1+c)(1+d) 
                \pmod{I^{i+1}_K(G)}.\\
\endsplit
$$
Case 2. If  $i\equiv 1\pmod{4}$,  then  we have 
$$
\split
{u}^{i}        &  \equiv(1+a)^{i}+
                \mu(1+a)^{i-2}(1+d) 
                \pmod{I^{i+1}_K(G)},\\
{v}{u}^{i-1}     &  \equiv(1+a)^{i-1}(1+c)+
                (1+a)^{i-2}(1+d) 
                \pmod{I^{i+1}_K(G)},\\
({u}{v}){u}^{i-2 } &  \equiv(1+a)^{i-1}(1+c) 
                \pmod{I^{i+1}_K(G)},\\
({v}{u}{v}){u}^{i-3} &  \equiv(1+a)^{i-3}(1+c)(1+d) 
                \pmod{I^{i+1}_K(G)}.\\
\endsplit
$$
Case 3. If  $i\equiv 2\pmod{4}$,  then  we have 
$$
\split
{u}^{i}        &  \equiv(1+a)^{i}
                  \pmod{I^{i+1}_K(G)},\\
{v}{u}^{i-1}   &  \equiv(1+a)^{i-1}(1+c)+
                  (1+a)^{i-2}(1+d)\\
               &  +\mu(1+a)^{i-3}(1+c)(1+d) 
                  \pmod{I^{i+1}_K(G)},\\
({u}{v}){u}^{i-2 } &  \equiv(1+a)^{i-1}(1+c)+
                  \mu(1+a)^{i-3}(1+c)(1+d) 
                  \pmod{I^{i+1}_K(G)},\\
({v}{u}{v}){u}^{i-3} &  \equiv(1+a)^{i-3}(1+c)(1+d) 
                \pmod{I^{i+1}_K(G)}.\\
\endsplit
$$
Case 4. If $i\equiv 3\pmod{4}$,  then we have
$$
\split
{u}^{i}          &  \equiv(1+a)^{i}+
                  \mu(1+a)^{i-2}(1+d) 
                  \pmod{I^{i+1}_K(G)},\\
{v}{u}^{i-1}       &  \equiv(1+a)^{i-1}(1+c)+
                  (1+a)^{i-2}(1+d) 
                  \pmod{I^{i+1}_K(G)},\\
({u}{v}){u}^{i-2 }   &  \equiv(1+a)^{i-1}(1+c)
                  \pmod{I^{i+1}_K(G)},\\
({v}{u}{v}){u}^{i-3}   &  \equiv(1+a)^{i-3}(1+c)(1+d) 
                  \pmod{I^{i+1}_K(G)}.\\
\endsplit
$$
It follows that  the  elements $\{\,\,\, {u}^{i},
{v}{u}^{i-1}, {u}{v}{u}^{i-2}, {v}{u}{v}{u}^{i-3}
\,\,\,\}$ are linearly independent  modulo ${I^{i+1}_K(G)}$.
Therefore, the matrix of decomposition is regular  and 
$\{1, u^i, vu^j, vuvu^k, uvu^l\mid 0\leq i,j,k,l \}$ form  a
filtered multiplicative $K$-basis of  $KG$.

Now let $G$ be  one of the  groups   $G_{13}$, $G_{14}$,
$G_{18}$, $G_{23}$, $G_{24}$ or  $G_{25}$ with $m>5$.  Clearly, $G$
is two-generated and we can assume that  $u=b_1$,  $v=b_2$, where
$b_1, b_2$  can be written as in  (4).   Then by (1) we get $1+d$ 
is central  modulo  $I^3_K(G)$ and 
$$
(1+c)(1+a)\equiv (1+a)(1+c)+(1+a)^2+(1+d)
                            \pmod{I^3_K(G)}.\tag5
$$
It follows  that 
$$
\cases
{u}{v} \equiv (\alpha_1+\alpha_2)\beta_1(1+a)^2+ 
\Delta(1+a)(1+c)+\alpha_2\beta_1(1+d)
&\pmod{I^3_K(G)};\\
{v}{u} \equiv \alpha_1(\beta_1+\beta_2)(1+a)^2+
\Delta(1+a)(1+c)+\alpha_1\beta_2(1+d)
&\pmod{I^3_K(G)};\\
{u}^2\equiv \alpha_1(\alpha_1+\alpha_2)(1+a)^2+\alpha_1\alpha_2(1+d)
&\pmod{I^3_K(G)};\\
{v}^2\equiv \beta_1(\beta_1+\beta_2) (1+a)^2+\beta_1 \beta_2 (1+d)
&\pmod{I^3_K(G)}.\\
\endcases
$$
Since the $K$-dimension of $I^2_K(G)/I^3_K(G)$ equals $3$ and
$\Delta\not= 0$, we have  that ${u}{v}\not\equiv
{v}{u}\pmod{I^3_K(G)}$ and ${u}^2\not\equiv
{v}^2\pmod{I^3_K(G)}$. Thus either ${v}^2\equiv 0
\pmod{I^3_K(G)}$ and $\beta_1=0$ or ${u}^2\equiv 0
\pmod{I^3_K(G)}$ and $\alpha_1=0$. It is easy to see  that
the second case is symmetric to the first one, so we consider
only the first one.  Therefore, ${u}= (1+a)+\mu(1+c)$ and
${v}= 1+c$, where $\mu=\frac{\alpha_2}{\alpha_1}$.

We will prove  that the elements $\{({u}{v})^{i}{u},
{u}^2{v}({u}{v})^{i-1}, ({v}{u})^{i}{v},
{u}^3({v}{u})^{i-1}\}$ form a basis  modulo $I^{2i+2}_K(G)$
and the elements $\{({u}{v})^{i}, {u}^2({v}{u})^{i-1},
({v}{u})^{i}, {u}^2({u}{v})^{i-1}\}$ form a basis   modulo
$I^{2i+1}_K(G)$, where $i>1$.
 
First, it is easy to see by  induction  that 
$$
\split
({u}{v})^{2i+1}       & \equiv (1+a)^{4i+1}(1+c) 
                        \pmod{I^{4i+3}_K(G)},\\
({u}{v})^{2i+2}       & \equiv (1+a)^{4i+3}(1+c)+ 
                        (1+a)^{4i+1}(1+c)(1+d)
                        \pmod{I^{4i+5}_K(G)},\\
({v}{u})^{2i+1}       & \equiv (1+a)^{4i+2}+
                        (1+a)^{4i+1}(1+c)\\
                      & +(1+a)^{4i}(1+d)  
                        \pmod{I^{4i+3}_K(G)},\\
({v}{u})^{2i+2}       & \equiv (1+a)^{4i+4}+
                        (1+a)^{4i+3}(1+c)\\
                      & +(1+a)^{4i+1}(1+c)(1+d) 
                        \pmod{I^{4i+5}_K(G)}. 
\endsplit
$$

By (1) and (5) we also have  modulo  ${I^{4}_K(G)}$ that 
$$
\split
(1+c)(1+a)^2  \equiv &(1+a)^2(1+c),\\
{u}^2         \equiv &(1+\mu)(1+a)^2+\mu(1+d),\\
{u}^2{v}      \equiv &
              (1+\mu)(1+a)^2(1+c)+
              \mu(1+c)(1+d),\\
{u}^3         \equiv &
              (1+\mu)(1+a)^3+\mu(1+\mu) (1+a)^2(1+c)\\
              &+ \mu(1+a)(1+d)+ \mu^2(1+c)(1+d).
\endsplit
$$
Now we consider  the following two cases:  

Case 1. Let $k=2j+1$ be  odd. Then 
$$
\split
({u}{v})^{k}{u} 
            \equiv ({u}{v})^{2j+1}{u} 
          & \equiv (1+a)^{4j+3} 
            +(1+a)^{4j+2}(1+c)\\ 
          & +(1+a)^{4j+1}(1+d) 
            \pmod{I^{4j+4}_K(G)},\\
 {u}^2{v}({u}{v})^{k-1} 
            \equiv {u}^2{v}({u}{v})^{2(j-1)+2}  
          & \equiv (1+\mu)(1+a)^{4j+2}(1+c)\\
          & +(1+a)^{4j}(1+c)(1+d) 
            \pmod{I^{4j+4}_K(G)},\\
\endsplit
$$$$
\split
 ({v}{u})^{k}{v} 
            \equiv ({v}{u})^{2j+1}{v} 
          & \equiv (1+a)^{4j+2}(1+c)\\
          & + (1+a)^{4j}(1+c)(1+d) 
            \pmod{I^{4j+4}_K(G)},\\
{u}^3({v}{u})^{k-1}
            \equiv {u}^3({v}{u})^{2(j-1)+2}  
          & \equiv (1+\mu)(1+a)^{4j+3}+
            (1+\mu)(1+a)^{4j+3}(1+c)\\
          & +(1+\mu)(1+a)^{4j}(1+c)(1+d)\\ 
          & +\mu(1+a)^{4j+1}(1+d)
            \pmod{I^{4j+4}_K(G)},\\
\endsplit
$$$$
\split          
 {u}^2({v}{u})^{k-1}
            \equiv {u}^2({v}{u})^{2(j-1)+2}   
          & \equiv (1+\mu)(1+a)^{4j+2}\\
          & +(1+\mu)(1+a)^{4j+1}(1+c)
            +\mu(1+a)^{4j}(1+d)\\
          & +(1+a)^{4j-1}(1+c)(1+d) 
            \pmod{I^{4j+3}_K(G)},\\
 ({v}{u})^{k}
            \equiv ({v}{u})^{2j+1} 
          & \equiv (1+a)^{4j+2}
            +(1+a)^{4j+1}(1+c)\\
          & + (1+a)^{4j}(1+d) 
            \pmod{I^{4j+3}_K(G)},\\
\endsplit
$$$$
\split            
 {u}^2({u}{v})^{k-1}
            \equiv {u}^2({u}{v})^{2(j-1)+2}  
          & \equiv (1+\mu)(1+a)^{4j+1}(1+c)\\
          & +(1+a)^{4j-1}(1+c)(1+d) 
            \pmod{I^{4j+3}_K(G)},\\
({u}{v})^{k} 
            \equiv ({u}{v})^{2j+1}   
          & \equiv (1+a)^{4j+1}(1+c) 
            \pmod{I^{4j+3}_K(G)}.\\
\endsplit
$$

Case 2. Let  $k=2j$ be  even. Then 
$$
\split
({u}{v})^{k}{u} 
            \equiv ({u}{v})^{2(j-1)+2}{u} 
          & \equiv (1+a)^{4j+1} 
            +(1+a)^{4j}(1+c)\\
          & +(1+a)^{4j-2}(1+c)(1+d) 
            \pmod{I^{4j+2}_K(G)},\\
{u}^2{v}({u}{v})^{k-1} 
            \equiv {u}^2{v}({u}{v})^{2(j-1)+1}
          &  \equiv (1+\mu)(1+a)^{4j}(1+c)\\
          & +(1+a)^{4j-2}(1+c)(1+d) 
            \pmod{I^{4j+2}_K(G)},\\
\endsplit
$$$$
\split
{u}^3({v}{u})^{k-1}
            \equiv {u}^3({v}{u})^{2(j-1)+1}  
          & \equiv (1+\mu)(1+a)^{4j+1}
            +(1+\mu)(1+a)^{4j}(1+c)\\ 
          & +\mu(1+\mu)(1+a)^{4j-2}(1+c)(1+d)\\ 
          & +(1+a)^{4j-1}(1+d)
            \pmod{I^{4j+2}_K(G)},\\
({v}{u})^{k}{v} 
            \equiv ({v}{u})^{2(j-1)+2}{v} 
          & \equiv (1+a)^{4j}(1+c) 
            \pmod{I^{4j+2}_K(G)},\\
\endsplit
$$$$
\split
({u}{v})^{k} 
            \equiv ({u}{v})^{2(j-1)+2}   
          & \equiv (1+a)^{4j-1}(1+c)\\
          & + (1+a)^{4j-3}(1+c)(1+d) 
            \pmod{I^{4j+1}_K(G)},\\
{u}^2({v}{u})^{k-1}
            \equiv {u}^2({v}{u})^{2(j-1)+1}   
          & \equiv (1+\mu)(1+a)^{4j}+
            (1+\mu)(1+a)^{4j-1}(1+c)\\
          & +(1+a)^{4j-2}(1+d)\\
          & +\mu(1+a)^{4j-3}(1+c)(1+d) 
            \pmod{I^{4j+1}_K(G)},\\
\endsplit
$$$$
\split
({v}{u})^{k}
            \equiv ({v}{u})^{2(j-1)+2}  
          & \equiv (1+a)^{4j} +
            (1+a)^{4j-1}(1+c)\\
          & +(1+a)^{4j-3}(1+c)(1+d) 
            \pmod{I^{4j+1}_K(G)},\\
{u}^2({u}{v})^{k-1}
            \equiv {u}^2({u}{v})^{2(j-1)+1}  
          & \equiv  (1+\mu)(1+a)^{4j-1}(1+c)\\
          & +\mu(1+a)^{4j-3}(1+c)(1+d) 
            \pmod{I^{4j+1}_K(G)}.\\
\endsplit
$$
It follows that $({u}{v})^{k}$, ${u}^2({v}{u})^{k-1}$,
$({v}{u})^{k}$, ${u}^2({u}{v})^{k-1}$ and also
$({v}{u})^{k}{u}$, ${u}^2v({u}{v})^{k-1}$, $({v}{u})^{k}v$,
${u}^3({v}{u})^{k-1}$ are linearly independent   modulo
$I^{2k+1}_K(G)$ and   modulo $I^{2k+1}_K(G)$, respectively.
Therefore, as before, the matrix of decomposition is regular
and
$$
\{ 
({u}{v})^{i}{u}, ({v}{u})^{i}{v}, ({u}{v})^{i}, ({v}{u})^{i}, 
{u}^2{v}({u}{v})^{j}, {u}^3({v}{u})^{j},  {u}^2({v}{u})^{j},
{u}^2({u}{v})^{j}
\}
$$  
form a filtered multiplicative $K$-basis of $KG$.

\enddocument